\magnification=\magstep1                          
\hsize=15truecm                                   
\hoffset=1.2cm                                    

\parskip=5truept plus 1truept minus 1truept       
\magnification \magstep1                          
\hsize=15truecm                                   
\hoffset=1.cm                                     

\parskip=5truept plus 1truept minus 1truept       
\baselineskip=20truept plus 1truept minus 1truept 

\font\grande=cmr10 scaled \magstep2               
\font\note=cmr8                                   
\font\maiuscoletto=cmcsc10 scaled \magstep1       

\def\titolo#1#2#3{{\centerline {\bf \grande {#1}}}\medskip
{\centerline {\bf \grande {#2}}}\medskip {\centerline {\bf \grande {#3}}}}
\def\abstract#1{\noindent {\bf Abstract.} {\narrower {\note #1}\par}}

\font\tenmsb=msbm10
\font\sevenmsb=msbm7
\font\fivemsb=msbm5
\newfam\msbfam
\textfont\msbfam=\tenmsb
\scriptfont\msbfam=\sevenmsb
\scriptscriptfont\msbfam=\fivemsb
\def\Bbb{\fam\msbfam }

\def\RR{{\Bbb R}}

\def\n{\noindent}
\def\vfi{\varphi}

\def\bull{{\vrule height.9ex width.8ex depth-.1ex}}       

\def\max{\mathop{\rm max}}
\def\sup{\mathop{\rm sup} }

\def\ve{\varepsilon}

\newdimen\shift \newbox\leftbox 
\def\newitem#1{\par\setbox\leftbox=\hbox{#1}\shift=\wd\leftbox
\def\rightshift{\hskip -\parindent\hskip\shift\hangindent\shift}
\everypar={\rightshift}
\leavevmode\hskip -\shift\box\leftbox}
\def\endnewitem{\everypar={}}                     
\newdimen\rientrodescr \rientrodescr=20pt

\pageno=-2
\nopagenumbers
\phantom{a} \vskip2truecm
\titolo {AN INTERIOR ESTIMATE} {FOR A NONLINEAR} {PARABOLIC EQUATION} 

\vskip1truecm                                      

\centerline {\maiuscoletto Giuseppe Maria
Coclite}
 \centerline  {\note S.I.S.S.A., via Beirut 2-4,
Trieste 34014, Italy.}
\centerline  {\note e-mail: coclite@sissa.it.}
\vskip2truecm 

\centerline  {Ref. S.I.S.S.A. 52/2002/M (July 2002)}

\vskip2truecm                                      

\abstract {In this paper there are estimated the 
derivatives of the solution of an initial boundary value problem for a
nonlinear uniformly parabolic equation in the interior with the total variation
of the boundary data and the $L^{\infty}-$norm of
the initial condition.}

\vfill
\eject

\quad\vfill\eject

\footline{\hfill\folio\hfill}
\pageno=1
\noindent {\bf 1. Introduction}
\bigskip

\bigskip
In this paper we are interested to estimate the $L^1-$norms of the derivatives
of the solution of a nonlinear parabolic equation. More precisely, we estimate
the derivatives of  the solution of an initial boundary value problem  in the
interior with respect to the total variation of the boundary data and the
$L^{\infty}-$norm of the initial datum. We split the solution in three parts,
one depending only on the initial datum and the other two depending only on
the boundary data. Moreover these maps are solutions of a linear parabolic
equation. The main tools of the proofs are the Maximum Principle and Energy
estimates.

Let $ W=W(t,x)$ the solution of the quasilinear initial boundary value problem (see [3, Theorem VI 5.2]) 

$$\cases{\displaystyle
 W_t = a (x, W) W_{xx} ,\quad & ${\rm for}\>\>  0\le x \le 1,\>\> t\ge 0 ,
$\cr  {}& {}\cr 
W(0,x)=\varphi (x), \quad & ${\rm for}\>\>0\le x \le 1, $\cr
{}& {}\cr
W (t,0)=g_0(t), \quad & ${\rm for}\>\> t\ge 0 ,$\cr
{}& {}\cr
W (t,1)=g_1(t), \quad & ${\rm for}\>\> t\ge 0 ,$\cr} \eqno(1)$$
where 
$$a \in C^3 (\RR^2), \quad 0 < a_* \le a (\cdot, \cdot) \le
a^* < +\infty, \quad \Vert a \Vert_{C^3} \le k \eqno(2)$$
and
$$\varphi \in C^2 ([0,1]), \quad g_0, g_1 \in  C^1 (\RR_+) \cap
BV(\RR_+),\quad g_0 (0)= \varphi (0), \>\>  g_1 (0)= \varphi
(1).\eqno(3)$$ The main results of this paper are  the following ones.
\bigskip

\bigskip
\noindent{\maiuscoletto Theorem 1. \phantom{a}}{\it Let} $ W= W( t,x)$
{\it be the classical solution of (1), with $a =a (x,y) $ satisfying
(2), $ \varphi= \varphi (x), \> g_0=g_0(t), \> g_1=g_1(t)$ satisfying (3) and 
$ c_1 >0 $. There
exists $ C=C(c_1, k, a_*, a^*,\Vert
\varphi\Vert_{L^{\infty}},\Vert g_0\Vert_{L^{\infty}},\Vert g_1
\Vert_{L^{\infty}})>0 $ such that} 
 $$\int_{c_1}^T 
\vert W_{t} (t,x)\vert dt \le C\>
\Vert \varphi \Vert_{L^{\infty}} + \int_{0}^T \Big( \vert g_0'(t)
\vert + \vert g_1'(t)\vert \Big) dt  \eqno(4)$$ 
{\it for all $0\le x\le 1$ and $ T \ge c_1 .$}
\bigskip

\noindent{\maiuscoletto Theorem 2. \phantom{a}}{\it Let} $ W= W(x, t)$
{\it be the classical solution of (1), with $a =a (x,y) $ satisfying
(2), $ \varphi= \varphi (x), \> g_0=g_0(t), \> g_1=g_1(t)$ satisfying (3),
$ c_1 >0 $ and $0 < \varepsilon < \displaystyle{1 \over 2}.$ There
exists $ C=C(\varepsilon, c_1, k, a_*, a^*,\Vert
\varphi\Vert_{L^{\infty}},\Vert g_0\Vert_{L^{\infty}},\Vert g_1
\Vert_{L^{\infty}})>0 $ such that} 
$$\int_{c_1}^T dt
\int_{\varepsilon}^{1-\varepsilon} \vert W_{tx} (t,x)\vert dx \le C
\Bigg(\Vert \varphi \Vert_{L^\infty} + \int_{0}^T \Big( \vert
g_0'(t) \vert + \vert g_1'(t)\vert \Big) dt \Bigg) \eqno(5)$$ 
{\it for all $ T \ge c_1 $.}

In the literature there are well-known interior estimates on the
$L^{\infty}-$norm of the derivatives of the solution of (1), called Schauder
estimates (e.g. see [2]) and the ones on the $L^2-$norm (e.g. see [4]).

As a motivation  and application  of this results we can see [1]. There the
authors prove the convergence of the vanishing viscosity solutions for a
particular $2\times2$ system of conservation laws. They show the compactness
of that family of solutions via uniform estimates on the total variation and
 Helly's Theorem. A basic ingredient  of these estimates (see [1, Lemma 3]) is
proved here as Theorem 2.

Denote
$$\bar a(t,x)\doteq a \big(x,W(t,x)\big),\qquad 0\le x\le 1,\>t\ge 0,$$
and consider the solution $u=u(t,x)$ of the problem
$$\cases{\displaystyle
 u_t = \bar a(t,x) u_{xx} ,\quad & ${\rm for}\>\>  0\le x \le 1,\>\> t\ge 0 ,
$\cr  {}& {}\cr 
u(0,x)=\varphi (x), \quad & ${\rm for}\>\>0\le x \le 1, $\cr
{}& {}\cr
u ( t,0)=g_0(t), \quad & ${\rm for}\>\> t\ge 0 ,$\cr
{}& {}\cr
u ( t,1)=g_1(t), \quad & ${\rm for}\>\> t\ge 0.$\cr} \eqno(6)$$
By uniqueness, $W\equiv u$ (see [3, Theorem VI 5.2]).  Now fix $ c_1 >0$ and $0 <
\varepsilon < \displaystyle{1 \over 2}$, by Schauder estimates (see [2,
Chapter 3, Section8]), there exists $K_0=K_0(c_1,\>\varepsilon,\>k, \> a_*,\> 
a^*)>0$ such that 
$$\sup\limits_{(t,\> x)\in S}\big\{\vert W_{x}\vert,\>\vert
W_{t}\vert,\>\vert W_{tx}\vert\big\}\le K_0 \big(\Vert \varphi\Vert_{L^{\infty}}
+\Vert g_0\Vert_{L^{\infty}}+\Vert g_1 \Vert_{L^{\infty}}\big),$$
and so, by the definition of $\bar a$ and (2), 
$$\sup\limits_{(t,\> x)\in S}\big\{\vert \bar a_{t}\vert,\>\vert \bar
a_{tx}\vert\big\}\le K_1 \big(\Vert \varphi\Vert_{L^{\infty}} +\Vert
g_0\Vert_{L^{\infty}}+\Vert g_1 \Vert_{L^{\infty}}\big),$$
for some constant $K_1=K_1(c_1,\>\varepsilon,\>k, \> a_*,\>  a^*)>0,$
where
$$S\doteq \big\{(t,x)\in \RR^2 ; c_1\le t,\>\> \varepsilon \le x \le
1-\varepsilon\big\}.$$

To simplify some technical aspects of the proofs,  we shall assume also that (see(3)) 
$$\varphi (0)= \varphi (1) =g_0 (0)=g_1 (0)=0.$$
Let $ u_1=u_1(t,x),\> u_2=u_2(t,x),\> u_3 =u_3(t,x)$ be
the solutions of the linear equation 
$$u_t = \bar a(t,x) u_{xx},\qquad t>0,\>\> 0<x<1, \eqno(7)$$
satisfying the initial and boundary conditions
$$u_1(0,\cdot) \equiv \varphi,\>\> u_1(\cdot, 0)\equiv0,\>\> u_1(\cdot,
1)\equiv0,$$
$$u_2(0,\cdot) \equiv 0,\>\> u_2(\cdot, 0)\equiv
g_0,\>\>u_2(\cdot, 1)\equiv0,$$
$$u_3(0,\cdot) \equiv 0,\>\> u_3(\cdot, 0)\equiv 0,\>\>u_3(\cdot,
1)\equiv g_1,$$
respectively. By (6) and the linearity of (7),
$$ W(t,x)=u(t,x) = u_1(t,x)+ u_2(t,x)+ u_3(t,x)\eqno(8)$$
is the solution of (1).

In Section 2 we prove the estimates (4) and (5) for $u_2$ and $u_3$, namely we
consider the case of (1) with the null initial condition. On the other side,
in Section 3 we prove the same ones for $u_1$, namely we consider (1) with the
null boundary data. In Section 4 we give the proofs of Theorems 1 and 2.
Finally in the last Section (the Appendix) we prove two lemmas, the first one
is a simply measure theory result, the second one consists of two Poincare type
inequalities. The proofs of these two lemmas are needed for the sake of the
best-constants.

\bigskip
\noindent {\bf 2. The Case with Null Initial Condition and General Boundary
Data} \bigskip

In this section we want to prove some estimates on the derivatives of the maps
$u_2$ and $u_3$, defined in the previous section.

\bigskip
\noindent{\maiuscoletto Lemma 3. \phantom{a}}{\it There results}
$$\int_0^T 
\vert u_{2,t} (t,x)\vert dt \le \int_0^T 
\vert g'_0 (t)\vert dt, $$ 
{\it for all  $0\le x \le 1.$}
\bigskip

\noindent{\maiuscoletto Proof.\phantom{a}} Define
$$\eqalign{h_0 (t) &\doteq {1\over 2 } \bigg( \int_0^t 
\vert g'_0 (\tau)\vert d\tau  + g_0 (t)\bigg), \cr k_0 (t) &\doteq {1\over 2
} \bigg( \int_0^t  \vert g'_0 (\tau)\vert d\tau  - g_0 (t)\bigg),}
\qquad t\ge 0.$$ 
Clearly $h_0,\, k_0 \in C^1(\RR_+) \cap BV(\RR_+),$  increasing in
$\RR_+$, positive and  
$$\eqalign{g_0 (t) &= h_0 (t)- k_0 (t), \cr\int_0^t 
\vert g'_0 (\tau)\vert d\tau &= \int_0^t 
 h'_0 (\tau) d\tau + \int_0^t 
 k'_0 (\tau) d\tau,} \qquad t\ge 0.\eqno(9)$$
Let $v_2$ and $\omega_2$ be the solutions of (7) such that
$$v_2(0,\cdot) \equiv 0,\>\> v_2(\cdot, 0)\equiv h_0,\>\>v_2(\cdot,
1)\equiv0,$$
$$\omega_2(0,\cdot) \equiv 0,\>\> \omega_2(\cdot, 0)\equiv
k_0,\>\>\omega_2(\cdot, 1)\equiv 0.$$
Since (7) is linear, by (9),
$$u_2 \equiv v_2 - \omega_2.\eqno(10)$$
Moreover, $v_{2,xx} $ and $\omega_{2,xx}$ are solutions of the equation
$$U_{t} = \bar a(t,x)U_{xx} + 2 \bar a_x(t,x) U_x + \bar a_{xx}(t,x)
U\eqno(11)$$ 
and, by the definition of $h_0$ and $k_0$,
$$v_{2,xx}(0,\cdot)\equiv v_{2,xx} (\cdot, 1) \equiv 0, \>\>v_{2,xx}(\cdot,
0)= {{h_0'}\over {\bar a(\cdot, 0)}}\ge 0, $$
$$\omega_{2,xx}(0,\cdot)\equiv \omega_{2,xx} (\cdot, 1) \equiv 0,
\>\>\omega_{2,xx}(\cdot, 0)= {{k_0'}\over {\bar a(\cdot, 0)}}\ge 0, $$
by the Maximum Principle (see [3, Theorem I
2.1]), $v_{2,xx} $ and $\omega_{2,xx}$ are positive. So $v_2 (t, \cdot)$ and
$\omega_2 (t,\cdot)$ are convex   in $[0,1],$ for each $t\ge 0.$
 By (2) and (7), we have
$$v_{2,t} (t,x) \ge 0, \quad \omega_{2,t} (t,x) \ge 0,\quad\quad t>0, \>\> 0\le x
\le 1.\eqno(12)$$
Fix $T\ge 0$ and $ 0\le x\le 1, $  by the Maximum Principle  and the
monotonicity of $h_0$, there results $$v_2 (T,x) \le \max\limits_{[0,T]} v_2
(\cdot, 0) = h_0 (T).$$ So, by (12) and since $h_0(0)=0$, we have 
$$\eqalign{\int_0^T \vert v_{2,t} (t,x)\vert dt &= \int_0^T  v_{2,t}
(t,x) dt =\cr &= v_2 (T,x) \le h_0 (T) = \int_0^T  h'_0(t) dt}\eqno(13)$$
and analogously
$$\int_0^T \vert \omega_{2,t} (t,x)\vert dt \le
\int_0^T  k'_0(t) dt ,\eqno(14)$$
then, by (9), (10), (13) and (14),
 $$\eqalign{ \int_0^T  \vert u_{2,t} (t,x)\vert dt &\le\int_0^T 
\big(\vert v_{2,t} (t,x)\vert + \vert \omega_{2,t} (t,x) \vert \big) dt\le \cr
& \le
\int_0^T  \big( h'_0(t)+ k'_0(t) \big) dt = \int_0^T 
\vert g'_0 (t)\vert dt.}$$
So the proof is concluded.\bull 

\bigskip
In the same way we can prove the following.
\bigskip

\noindent{\maiuscoletto Lemma 4. \phantom{a}}{\it There results}
$$\int_0^T 
\vert u_{3,t} (t,x)\vert dt \le \int_0^T 
\vert g'_1 (t)\vert dt, $$ 
{\it for all  $0\le x \le 1.$}
\bigskip

\noindent{\maiuscoletto Lemma 5. \phantom{a}}  {\it There exist two constants
$C_1,\> \delta_1>0$ depending only on $ a_*,$ $a^*,$ $k,$ 
$\ve,$  $c_1,$ $\Vert\varphi\Vert_{L^{\infty}},$ $\Vert
g_0\Vert_{L^{\infty}},$ $\Vert g_1 \Vert_{L^{\infty}}$
such that, if $\ve \le x_1<x_2\le 1-\ve$ and $x_2
-x_1 <\delta_1$, there results}
$$\int_{c_1}^T dt \int_{x_1}^{x_2} \vert
 u_{i,xt} ( t,x)\vert dx \le C_1 \int_{0}^T \Big( \vert u_{i,t} (t,x_1)
\vert + \vert u_{i,t} (t,x_2)\vert \Big) dt, $$ 
{\it for $i=2,3$ and all $ T \ge c_1.$}
\bigskip

\noindent{\maiuscoletto Proof.\phantom{a}} Fix $i\in\{2,\>3\}$ and $ 0\le t \le
T,$ by (7), there results 
$$\eqalign{&{d \over {dt}} \int_{x_1}^{x_2} {{ u_{i,tx}^2
(t,x)} \over 2} dx =  \int_{x_1}^{x_2}  u_{i,tx} (t,x)  u_{i,ttx} (t,x) dx =\cr
&=  u_{i,tx}(t,
x_2) u_{i,tt}(t, x_2) -   u_{i,tx}(t, x_1) u_{i,tt}(t, x_1) -\cr
&- \int_{x_1}^{x_2}  u_{i,txx} (t,x)  u_{i,tt} (t,x) dx =\cr
&=  u_{i,tx}(t,
x_2) u_{i,tt}(t, x_2) -  u_{i,tx}(t, x_1) u_{i,tt}(t, x_1) -\cr
&- \int_{x_1}^{x_2}  u_{i,txx} (t,x) \Big(\bar a_t (t,x)  u_{i,xx}
(t,x) + \bar a(t,x)  u_{i,txx} (t,x) \Big) dx=  \cr
&=  u_{i,tx}(t, x_2) u_{i,tt}(t, x_2) -   u_{i,tx}(t, x_1) u_{i,tt}(t,
x_1) -\cr
&- \int_{x_1}^{x_2} \bar a_t (t,x) u_{i,txx} (t,x) u_{i,xx} (t,x) dx 
- \int_{x_1}^{x_2} \bar a(t,x)  u_{i,txx}^2 (t,x)  dx\le\cr 
&\le u_{i,tx}(t, x_2) u_{i,tt}(t, x_2) -   u_{i,tx}(t, x_1)
u_{i,tt}(t, x_1) - \cr 
&-\int_{x_1}^{x_2} \bar a_t (t,x)u_{i,txx}
(t,x)  u_{i,xx} (t,x) dx- a_* \int_{x_1}^{x_2}   u_{i,txx}^2 (t,x)
 dx\le \cr
&\le  u_{i,tx}(t, x_2) u_{i,tt}(t, x_2) -   u_{i,tx}(t, x_1) u_{i,tt}(t,
x_1) +\cr
&+{1\over {2a_*}} \int_{x_1}^{x_2} \bar a_t^2 (t,x)   u_{i,xx}^2 (t,x) dx 
- {{a_*}\over 2} \int_{x_1}^{x_2} u_{i,txx}^2 (t,x)dx\le
\cr
&\le  u_{i,tx}(t, x_2) u_{i,tt}(t, x_2) -   u_{i,tx}(t, x_1) u_{i,tt}(t,
x_1) +\cr
&+{1\over {2a_*}} \int_{x_1}^{x_2} {{\bar a_t^2 (t,x)}\over {\bar a^2(t,x)}}  
u_{i,t}^2 (t,x) dx  - {{a_*}\over 2}
\int_{x_1}^{x_2} u_{i,txx}^2 (t,x)dx\le \cr
&\le  u_{i,tx}(t, x_2) u_{i,tt}(t, x_2) -   u_{i,tx}(t, x_1) u_{i,tt}(t,
x_1) +\cr
&+{{\Vert \bar a_t\Vert_{L^\infty(S) }^{2}}\over {2a_*^3}} \int_{x_1}^{x_2}
  u_{i,t}^2 (t,x) dx  - {{a_*}\over 2}
\int_{x_1}^{x_2} u_{i,txx}^2 (t,x)dx.}$$ 

\n By formulas (A.2) and (A.3) of the Appendix, we have
$$-\int_{x_1}^{x_2}  u_{i,txx}^2 (t,x)  dx \le - \int_{x_1}^{x_2}  {u_{i,tx}^2 (t,x)\over
{2(x_2-x_1)^2}}  dx  + 
{{\vert u_{i,t} (t,x_2) -u_{i,t} (t,x_1) \vert^2} \over {|x_2 -x_1|^3}},$$
$$\int_{x_1}^{x_2}  u_{i,t}^2 (t,x)  dx \le 2(x_2-x_1)^2\int_{x_1}^{x_2} 
u_{i,tx}^2 (t,x) dx+ 2(x_2 - x_1) \vert u_{i,t}(t,x_1)\vert^2 ,
$$
respectively and then
$$\eqalign{&{d \over {dt}} \int_{x_1}^{x_2} {{u_{i,tx}^2 (t,x)} \over {2
}} dx \le\cr
&\le u_{i,tx}(t, x_2)u_{i,tt}(t, x_2) -  u_{i,tx}(t, x_1)u_{i,tt}(t,
x_1)+\cr
&+\left[-{{a_*}\over{4 (x_2-x_1)^2}} + {{\Vert \bar a_t\Vert_{L^\infty(S)
}^{2}}\over {a_*^3}}  (x_2-x_1)^2\right]\int_{x_1}^{x_2}  u_{i,tx}^2 (t,x) 
dx+\cr &+{{\Vert \bar a_t\Vert_{L^\infty(S) }^{2}}\over {a_*^3}} (x_2-x_1)\vert
u_{i,t}(t,x_1)\vert^2 + {{a_*}\over 2} {{\vert u_{i,t} (t,x_2) -u_{i,t} (t,x_1)
\vert^2} \over {|x_2 -x_1|^3}}.}\eqno(15)$$
Moreover, there exists  $\delta_1>0$ such that, if $x_2
-x_1 <\delta_1$,
there results
$$\eqalign{& u_{i,tx}(t, x_2)u_{i,tt}(t, x_2) -  u_{i,tx}(t, x_1)u_{i,tt}(t,
x_1)+\cr
&+ {{\Vert \bar a_t\Vert_{L^\infty(S) }^{2}}\over {a_*^3}} (x_2-x_1)^2 \int_{x_1}^{x_2}  u_{i,tx}^2 (t,x) 
dx+\cr
&+{{\Vert \bar a_t\Vert_{L^\infty(S) }^{2}}\over {a_*^3}} (x_2-x_1)\vert
u_{i,t}(t,x_1)\vert^2 \le\cr 
&\le {{a_*}\over{8 (x_2-x_1)^2}}\int_{x_1}^{x_2} 
u_{i,tx}^2 (t,x)  dx +{a_* \over 2}
{{\vert u_{i,t} (t,x_2) -u_{i,t} (t,x_1) \vert^2} \over {|x_2 -x_1|^3}},}
$$
and then, by (15),
$$\eqalign{&{d \over {dt}} \int_{x_1}^{x_2} {{u_{i,tx}^2 (t,x)}
\over 2} dx \le\cr
&\le - {{a_*}\over 8}\int_{x_1}^{x_2} 
{{u_{i,tx}^2 (t,x) }\over { (x_2-x_1)^2}}  dx  + a_* {{\vert u_{i,t}
(t,x_2) -u_{i,t} (t,x_1) \vert^2} \over {|x_2 -x_1|^3}},}$$
 hence 
$$\eqalign{&\int_{x_1}^{x_2} {u_{i,tx}^2 (t,x)} dx \le\cr
&\le 2 a_*\displaystyle{{e^{-{{   a_* t} \over {4(x_2 -x_1)^2}}}}\over {|x_2
-x_1|^3}} \int_0^t e^{{{  a_* \tau} \over {4(x_2 -x_1)^2}}}  \vert u_{i,t}
(\tau,x_2) - u_{i,t} (\tau,x_1) \vert^2 d\tau.}
\eqno(16)$$
Since 
$$\Bigg(\int_{x_1}^{x_2} \vert u_{i,tx} (t,x)\vert dx\Bigg)^2
\le(x_2-x_1)\int_{x_1}^{x_2} {u_{i,tx}^2 (t,x)} dx,\eqno(17)$$
and
$$  \vert u_{i,t} (\tau,x_2)
- u_{i,t} (\tau,x_1) \vert \le\int_{x_1}^{x_2} \vert u_{i,tx}
(\tau,x)\vert dx,\eqno(18)$$
we have, by (16), (17) and (18),
$$ f^2(t) \le e^{-\lambda t} \int_0^t e^{\lambda \tau} h(\tau)
f(\tau)d\tau,$$
where
$$\eqalign{f(t)&\doteq \int_{x_1}^{x_2} \vert u_{i,tx} (t,x)\vert dx,\cr
 h(t)&\doteq 2a_*\displaystyle{{\vert u_{i,t}
(t,x_2) - u_{i,t} (t,x_1) \vert} \over {|x_2 -x_1|^2}},\cr
 \lambda &\doteq
{{   a_*} \over {4(x_2 -x_1)^2}}.}$$  
By Lemma A-1 of the Appendix, there results 
$$\int_0^Tdt\int_{x_1}^{x_2} \vert u_{i,tx} (t,x)\vert  dx  \le
C_1\int_0^T \vert u_{i,t} (\tau,x_2) - u_{i,t}
(\tau,x_1) \vert dt.$$ 
So the proof is concluded.\bull

\bigskip
\noindent {\bf 3. The Case with Homogeneous Dirichlet Boundary Conditions}  
\bigskip

In this section we want to prove some estimates on the derivatives of $u_1$,
defined in Section 1.

\bigskip
\noindent{\maiuscoletto Lemma 6. \phantom{a}} {\it There exists a constants
$C_2\> >0$ depending only on $ a_*,$ $a^*,$ $k,$  $c_1,$ $\Vert\varphi\Vert_{L^{\infty}},$ $\Vert
g_0\Vert_{L^{\infty}},$ $\Vert g_1 \Vert_{L^{\infty}}$
such that, there results}
$$\int_{c_1}^T 
\vert u_{1,t} (t,x)\vert dt \le C_2 \Vert\varphi
\Vert_{L^\infty}, $$ 
{\it for all $0<c_1 \le T$ and $0\le x \le 1.$}
\bigskip

\noindent{\maiuscoletto Proof.\phantom{a}}  Let $v_1=v_1(t,x)$ and $\omega_1=\omega_1(t,x)$ be the solutions of (7) 
satisfying the following conditions
$$\eqalign{v_1(c_1,x)&= {1\over 2} \Big( \int_0^x\int_0^y | u_{1,xx}(c_1,\xi)|d\xi dy +  u_{1}(c_1,x)
\Big),\cr
\omega_1(c_1,x) &= {1\over 2}\Big( \int_0^x \int_0^y | u_{1,xx}(c_1,\xi)|d\xi dy - u_{1}(c_1,x)
\Big),}\eqno(19)$$
for $0\le x\le 1$ and 
$$v_1(t,0)= v_1(t,1)= \omega_1(t,0)=\omega_1(t,1)=0,\eqno(20)$$
for   $t\ge c_1$.
By the linearity of (7), (19) and (20) there results
$$ u_1 = v_1-\omega_1,\qquad 0\le x\le 1,\>\>t\ge c_1.\eqno(21)$$
Moreover $v_1(c_1,\cdot)$ and $\omega_1(c_1,\cdot)$ are convex. Since $v_{1,xx}$ and $\omega_{1,xx}$ are solutions of (11)
and 
$$v_{1,xx}(c_1,\cdot)\ge 0,\quad v_{1,xx}(\cdot,0)={{v_{1,t}(\cdot,0)}\over {\bar a(\cdot,0)}}= 0,\quad
v_{1,xx}(\cdot,1)={{v_{1,t}(\cdot,1)}\over {\bar a(\cdot,1)}}= 0,$$
$$\omega_{1,xx}(c_1,\cdot)\ge 0,\quad \omega_{1,xx}(\cdot,0)={{\omega_{1,t}(\cdot,0)}\over {\bar a(\cdot,0)}}= 0,\quad
\omega_{1,xx}(\cdot,1)={{\omega_{1,t}(\cdot,1)}\over {\bar a(\cdot,1)}}= 0,$$
by the Maximum Principle, $v_{1,xx}$ and $\omega_{1,xx}$ are positive. So $v_{1}(t,\cdot)$ and $\omega_{1}(t,\cdot)$
 are convex in $[0,1]$,
for each $t\ge c_1$. By (2), (7), (19) and (20), we have 
$$v_{1,t}(t,x)\ge 0,\qquad \omega_{1,t}(t,x)\ge0,\qquad 0\le x\le 1,\>\>t\ge c_1.$$
Fix $T\ge c_1$ and $0\le x \le 1$, there results
$$\eqalign{\int_{c_1}^T|v_{1,t}(t,x)|dt &= \int_{c_1}^Tv_{1,t}(t,x)dt =\cr
&= v_1(T,x) -v_1(c_1,x) \le 2\Vert v_1\Vert_{L^{\infty}(S)}}\eqno(22)$$
and analogously
$$\int_{c_1}^T|\omega_{1,t}(t,x)|dt \le 2\Vert \omega_1\Vert_{L^{\infty}(S)}.\eqno(23)$$
By the Maximum Principle and the definition of $u_1, \>v_1,\>\omega_1$ there results
$$\eqalign{\Vert v_1\Vert_{L^{\infty}(S)} &+ \Vert \omega_1\Vert_{L^{\infty}(S)}\le\cr
&\le \Vert u_{1,xx}(c_1,\cdot)\Vert_{L^{\infty}([0,1])} + 
 \Vert u_{1}(c_1,\cdot)\Vert_{L^{\infty}([0,1])} .}$$
By the Maximum Principle and the definition of $u_1$, there results
$$ \Vert u_{1}(c_1,\cdot)\Vert_{L^{\infty}([0,1])}\le  \Vert \vfi\Vert_{L^{\infty}}.$$
Since 
$$u_{1,xx}(c_1,0)= {{u_{1,t}(c_1,0)}\over {\bar a(c_1,0)}}= 0,\quad u_{1,xx}(c_1,1)= {{u_{1,t}(c_1,1)}\over {\bar a(c_1,1)}}= 0,$$
there exists $0< \bar x <1$, depending on $c_1$, such that
$$ \Vert u_{1,xx}(c_1,\cdot)\Vert_{L^{\infty}([0,1])} = |u_{1,xx}(c_1,\bar x)|.$$
Let be $0<\bar \ve <\displaystyle{1\over 2}$ such that
$\bar \ve  \le \bar x \le 1-\bar \ve$, there results
$$ \Vert u_{1,xx}(c_1,\cdot)\Vert_{L^{\infty}([0,1])} = |u_{1,xx}(c_1,\bar x)| = 
\Vert u_{1,xx}(c_1,\cdot)\Vert_{L^{\infty}([\bar \ve, 1-\bar\ve)])}.$$
Moreover, by the Schauder Estimates, there exists a constant $K_2>0$ such that
$$ \Vert u_{1,xx}(c_1,\cdot)\Vert_{L^{\infty}([\bar \ve, 1-\bar\ve)])}\le K_2 \Vert \vfi \Vert_{L^{\infty}}
.\eqno(24)$$
Finally, by (21), (22), (23) and (24)  we can conclude
$$\eqalign{\int_{c_1}^T |u_{1,t}(t,x)|dt
&\le \int_{c_1}^T |v_{1,t}(t,x)|dt+\int_{c_1}^T |\omega_{1,t}(t,x)|dt\le\cr
&\le 2\big(\Vert v_1\Vert_{L^{\infty}(S)} + \Vert \omega_1\Vert_{L^{\infty}(S)}\big)\le 2 (  K_2 +1) \Vert \vfi\Vert_{L^{\infty}}.}$$
Since $K_2$ depends on $\bar\ve$ that depends on $c_1$,  the proof is done.\bull 

\bigskip

\noindent{\maiuscoletto Lemma 7. \phantom{a}}  {\it There exist two constants
$C_3, \delta_2 >0$ depending only on $ a_*,$ $a^*,$ $k,$ 
$\ve,$  $c_1,$ $\Vert\varphi\Vert_{L^{\infty}},$ $\Vert
g_0\Vert_{L^{\infty}},$ $\Vert g_1 \Vert_{L^{\infty}}$
such that, if $\ve \le x_1<x_2\le 1-\ve$ and  $x_2
-x_1 <\delta_2$, there results}
$$\int_{c_1}^T dt \int_{x_1}^{x_2} \vert
 u_{1,tx} (x, t)\vert dx \le C_3 \Vert \varphi\Vert_{L^\infty},$$  {\it for
all $ T \ge c_1.$} \bigskip

\noindent{\maiuscoletto Proof.\phantom{a}}
Call
$$\delta_2\doteq {{a_*}\over {8^{1/4}\Vert \bar a _t
\Vert_{L^{\infty}(S)}^{1/2}}},\eqno(25)$$ fix $\ve \le x_1 < x_2 \le 1-\ve$ and
consider the restriction of $u_1$ to the the strip 
$$\widetilde S \doteq\{(t,x)\in \RR^2;\> t\ge 0, \>\>x_1\le x \le x_2
\}.$$ There results  
$$u_1\equiv  \bar u +\widetilde u,\qquad {\rm in}\>\> \widetilde S,\eqno(26)$$
where $\bar u,\> \tilde u$ are the solutions of (7) such that 
$$\bar u(0,\cdot) \equiv
\varphi\big\vert_{[x_1,x_2]},\>\> \bar u(\cdot, x_1)\equiv\vfi(x_1),\>\> \bar u(\cdot,
x_2)\equiv \vfi(x_2),\eqno(27)$$
$$\tilde u(0,\cdot) \equiv 0,\>\> \tilde u(\cdot,
x_1)\equiv u_1(\cdot, x_1) -\vfi(x_1),\>\>\tilde u(\cdot, x_2)\equiv u_1(\cdot,
x_2)-\vfi(x_2).\eqno(28)$$
 By Lemmas 5 and 6, there results
$$\eqalign{&\int_{c_1}^T dt \int_{x_1}^{x_2} \vert
 \widetilde u_{tx} (t,x)\vert dx \le\cr &\le C_1\int_0^T \big(|\widetilde u_{t}
(t,x_1)|+|\widetilde u_{t} (t,x_2)|\big)dt=\cr
&=2 C_1 \int_0^T \big(|u_{1,t} (t,x_1)| + |u_{1,t} (t,x_2)|\big)
\le 4 C_1C_2 \Vert \vfi
\Vert_{L^\infty}.}\eqno(29)$$
Finally, denote
$$\hat u(t,x) \doteq \bar u(t,x) -{{x-x_1}\over{x_2-x_1}}\vfi(x_2) -{{x_2-x}\over{x_2-x_1}}\vfi(x_1). $$
Clearly $\hat u$ is solution of (7) and there results
$$\hat u(0,x)= \vfi(x) -{{x-x_1}\over{x_2-x_1}}\vfi(x_2) -{{x_2-x}\over{x_2-x_1}}\vfi(x_1),\eqno(30)$$
$$\hat u(t,x_1)= \hat u(t,x_2)=0,\eqno(31)$$
$$\hat u_{tx}(t,x)= \bar u_{tx}(t,x).\eqno(32) $$

Moreover, by the definition of $\hat u$, we obtain 
$$\eqalign{&{d \over {dt}} \int_{x_1}^{x_2} {{ \hat u_{tx}^2
(t,x)} \over 2} dx =  \int_{x_1}^{x_2}  \hat u_{tx} (t,x) \hat  u_{ttx} (t,x) dx =\cr
&=- \int_{x_1}^{x_2}  \hat u_{txx} (t,x) \hat u_{tt} (t,x) dx =\cr
&=  - \int_{x_1}^{x_2} \hat u_{txx} (t,x) \Big(\bar a_t (t,x) \hat u_{xx}
(t,x) + \bar a(t,x) \hat u_{txx} (t,x) \Big) dx=  \cr
&= - \int_{x_1}^{x_2} \bar a_t (t,x) \hat u_{txx} (t,x) \hat u_{xx} (t,x) dx
- \int_{x_1}^{x_2} \bar a(t,x) \hat u_{txx}^2 (t,x)  dx\le\cr 
&\le -\int_{x_1}^{x_2} \bar a_t (t,x)\hat u_{txx}
(t,x)\hat   u_{xx} (t,x) dx- a_* \int_{x_1}^{x_2}  \hat  u_{txx}^2 (t,x)
 dx\le \cr
&\le {1\over {2a_*}} \int_{x_1}^{x_2} \bar a_t^2 (t,x)  \hat u_{xx}^2 (t,x)
dx  - {{a_*}\over 2} \int_{x_1}^{x_2}\hat u_{txx}^2 (t,x)dx\le
\cr
&\le {1\over {2a_*}} \int_{x_1}^{x_2} {{\bar a_t^2 (t,x)}\over {\bar a^2(t,x)}}  
\hat u_{t}^2 (t,x) dx  - {{a_*}\over 2}
\int_{x_1}^{x_2} \hat u_{txx}^2 (t,x)dx\le \cr
&\le{{\Vert \bar a_t\Vert_{L^\infty(S) }^{2}}\over {2a_*^3}} \int_{x_1}^{x_2}
\hat  u_{t}^2 (t,x) dx  - {{a_*}\over 2}
\int_{x_1}^{x_2} \hat u_{txx}^2 (t,x)dx.}$$ 
 By (A.3) and (30), we have
$$-\int_{x_1}^{x_2}  \hat u_{txx}^2 (t,x)  dx \le - \int_{x_1}^{x_2} 
{\hat u_{tx}^2 (t,x)\over {2(x_2-x_1)^2}}  dx  ,$$
$$\int_{x_1}^{x_2} \hat  u_{t}^2 (t,x)  dx \le 2(x_2-x_1)^2\int_{x_1}^{x_2} 
\hat u_{tx}^2 (t,x),$$
 and then, by (25), if $x_2-x_1 \le \delta_2$,
$$\eqalign{&{d \over {dt}} \int_{x_1}^{x_2} {{\hat u_{tx}^2 (t,x)} \over {2
}} dx \le\cr
&\le\left[-{{a_*}\over{4 (x_2-x_1)^2}} + {{\Vert \bar a_t\Vert_{L^\infty(S)
}^{2}}\over {a_*^3}}  (x_2-x_1)^2\right]\int_{x_1}^{x_2}  \hat u_{tx}^2 (t,x) 
dx\le\cr
&\le\left[-{{a_*}\over{4 \delta_2^2}} + {{\Vert \bar a_t\Vert_{L^\infty(S)
}^{2}}\over {a_*^3}}  \delta_2^2\right]\int_{x_1}^{x_2}  \hat u_{tx}^2
(t,x)  dx\le\cr
&\le-{{a_*}\over{8 \delta_2^2}}\int_{x_1}^{x_2}  \hat u_{tx}^2 (t,x) 
dx.}$$  
Hence, by Schauder estimates and (30), there exists a constant $K_3>0$ such that
$$\eqalign{&\int_{x_1}^{x_2} {\hat u_{tx}^2 (t,x)} dx \le\cr
&\le \displaystyle{e^{-{{   a_* (t-c_1)}  \over {4\delta_2^2}}}}
 \int_{x_1}^{x_2} {\hat u_{tx}^2 (c_1,x)} dx\le\cr
&\le K_3\displaystyle e^{-{{   a_* t} \over {4\delta_2^2}}}\Vert \hat u(0,\cdot)
\Vert_{L^{\infty}}^2\le 9K_3\displaystyle e^{-{{   a_* t} \over {4\delta_2^2}}}\Vert\vfi
\Vert_{L^{\infty}}^2 ,}$$  
so 
$$\eqalign{&\int_{x_1}^{x_2} \vert \hat u_{tx} (t,x)\vert dx
\le(x_2-x_1)^{1\over 2}\left(\int_{x_1}^{x_2} {\hat u_{tx}^2 (t,x)}
dx\right)^{1\over 2} \le\cr
&\le 3 {K_3}^{1\over 2} \Vert \vfi
\Vert_{L^{\infty}}(x_2-x_1)^{1\over 2}\cdot e^{-{{   a_* t} \over
{8\delta_2^2}}}\le 3 {K_3}^{1\over 2}\Vert \vfi
\Vert_{L^{\infty}}\delta_2^{1\over 2}\cdot e^{-{{   a_* t} \over
{8\delta_2^2}}},}$$
 and integrating  on $[c_1,\> T]$
$$\int_{c_1}^T dt \int_{x_1}^{x_2} \vert \hat u_{tx} (t,x)\vert dx
\le{{24{K_3}^{1\over 2}\Vert \vfi
\Vert_{L^{\infty}}\delta_2^{5\over 2} }\over {a_*}}. \eqno(33)$$
By (26), we have
$$\int_{c_1}^T dt \int_{x_1}^{x_2} \vert u_{1,tx} (t,x)\vert dx
\le\int_{c_1}^T dt \int_{x_1}^{x_2} \big(\vert \bar u_{tx} (t,x)\vert  +
 \vert \tilde u_{tx} (t,x)\vert\big) dx,$$
then, by (29), (32) and (33), the thesis is done.\bull

\bigskip
\noindent {\bf 4. Proofs of Theorems 1 and 2}  
\bigskip

In this section we give the proofs of the main results of the paper.
 
\noindent{\maiuscoletto Proof of Theorem 1.\phantom{a}} The thesis
is direct consequence of  (8) and Lemmas 3, 4 and 6.\bull
\bigskip

\noindent{\maiuscoletto Proof of Theorem 2.\phantom{a}}
 Fix $ 0<c_1 \le
T, \> 0 <\varepsilon<\displaystyle{1\over 2}$ and observe that 
$$\eqalign{\int_{c_1}^T dt \int_{\varepsilon}^{1-\varepsilon} 
 \vert W_{tx} (t, x)\vert dx&=\int_{c_1}^T dt
\int_{\varepsilon}^{1-\varepsilon}   \vert u_{tx} ( t,x)\vert dx\le
\cr &\le\sum_{i=1,2,3}\int_{c_1}^T dt \int_{\varepsilon}^{1-\varepsilon}   \vert
u_{i,tx} ( t,x)\vert dx,}$$
where $W=W(t,x)$ and $u=u(t,x)$ are the solutions of (1) and (6), respectively.

Let $x_0,..., x_h$ such that
$$\varepsilon = x_0 < x_1 <....< x_{h-1}< x_h = 1-\varepsilon$$
and 
$$x_j -x_{j-1} <\delta, \qquad j=1,..., h,$$
where $\delta\le \min\{\delta_1,\> \delta_2\}$ and $\delta_1,\> \delta_2$ are
the ones of Lemmas 5 and 7, respectively. 

\n By Lemma 7, 
$$\eqalign{\int_{c_1}^T dt
\int_{\varepsilon}^{1-\varepsilon}  \vert u_{1,tx} (t,x)\vert dx&\le
\sum\limits_{j=1}^h\int_{c_1}^T dt\int_{x_{j-1}}^{x_j}   \vert u_{1,tx} (t,x)
\vert dx\le \cr &\le\sum\limits_{j=1}^h C_3\Vert \varphi \Vert_{L^\infty}= h 
C_3\Vert \varphi \Vert_{L^\infty}.} \eqno(34)$$
By Lemmas 3, 4 and 5, 
$$ \eqalign{&\sum_{i=2,3}\int_{0}^T dt
\int_{\varepsilon}^{1-\varepsilon}  \vert u_{i,tx} (t,x)\vert
dx=\cr
&=\sum_{i=2,3} \sum\limits_{j=1}^h \int_{0}^T dt
\int_{x_{j-1}}^{x_j}   \vert u_{i,tx} (t,x)\vert dx \le\cr
&\le \sum_{i=2,3}\sum\limits_{j=1}^h C_1 \int_{0}^T \Big( \vert
u_{i,t} (t,x_{j-1}) \vert + \vert u_{i,t} (t,x_j)\vert \Big) dt\le\cr
&\le 2\, h \cdot  C_1 \int_{0}^T \Big( \vert u_{t}
(t,0) \vert + \vert u_{t} (t,1)\vert \Big) dt= \cr
&=2 \, h \cdot  C_1
\int_{0}^T \Big( \vert g_0'(t) \vert + \vert g_1'(t)\vert \Big)
dt.}\eqno(35)$$
Since $h$, $C_1$ and $C_3$ depend only on $ a_*,$ $a^*,$ $k,$  
$\ve,$  $c_1,$ $\Vert\varphi\Vert_{L^{\infty}},$ $\Vert
g_0\Vert_{L^{\infty}},$ $\Vert g_1 \Vert_{L^{\infty}}$ the thesis is direct
consequence of  (34) and (35).\bull  \bigskip

\bigskip
\noindent {\bf Appendix: Two Technical Lemmas}  
\bigskip

\bigskip
In these section we  prove  two  lemmas.
The proofs of these two lemmas are more or less well-known. We insert here one
of  these ones for the sake of completeness and the best-constants.

\bigskip

\bigskip

\noindent{\maiuscoletto Lemma A-1. \phantom{a}}  {\it Let $f,\>h \,\in\>
C(\RR)$ be positive functions and fix a constant $\lambda >0$. If}
$$ f^2(t) \le e^{-\lambda t} \int_0^t e^{\lambda \tau} h(\tau)
f(\tau)d\tau,\qquad t\ge 0, \eqno(A.1)$$
{\it then}
$$\int_0^T f(t) dt \le {2\over \lambda} \int_0^T h(t) dt,$$
{\it for each $T\ge0.$}
\bigskip

\noindent{\maiuscoletto Lemma A-2. (Poincare type Inequalities )\phantom{a}} 
{\it For any $f\in C^2(\RR)$ one has}
$$\eqalignno{\int_a^b {f'}^2 (x) \,dx &\le 2\, (b-a)^2 \int_a^b {f''}^2 (x)
\,dx +2\, {{\big(f(b)-f(a) \big)^2} \over {b-a}},&\hbox{($A.2$)}\cr
\int_a^b {f}^2 (x) \,dx &\le 2\, (b-a)^2 \int_a^b {f'}^2 (x) \,dx +2\,
(b-a)\vert f(a)\vert^2,&\hbox{($A.3$)}\cr}$$  
%
{\it for each $ -\infty< a<
b<+\infty.$} \bigskip

\bigskip

\noindent{\maiuscoletto Proof of Lemma A-1.\phantom{a}} Fix $T\ge 0$. Define
$$g(t) \doteq e^{{\lambda \over 2} t} f(t), \qquad t\ge 0,$$
by (A.1), we have
$$ g^2(t) \le  \int_0^t e^{{\lambda \over 2} \tau} h(\tau)
g(\tau)d\tau,\qquad t\ge 0. \eqno(A.4)$$
Denote
$$M_t \doteq \sup\limits_{0\le \tau \le t} g(\tau), \qquad t\ge 0,$$
we claim that
$$M_t \le \int_0^t e^{{\lambda \over 2} \tau} h(\tau) d\tau,\qquad t\ge 0.
\eqno(A.5)$$
Fix $t \ge 0$. If $M_t=0$ we are done, assume that $M_t>0.$
Since $g$ is continuous there exists $0\le t_0\le t$ such that
$$M_t =g(t_0),$$
by (A.4) and the positivity of $f$ and $h$, there results
$$\eqalign{M_t^2 = g^2 (t_0)&\le \int_0^{t_0} e^{{\lambda \over 2} \tau} h(\tau)
g(\tau)d\tau \le\cr
&\le \int_0^{t} e^{{\lambda \over 2} \tau} h(\tau)
g(\tau)d\tau \le M_t \int_0^{t} e^{{\lambda \over 2} \tau} h(\tau)d\tau,}$$
so (A.5) is proved. Since
$$e^{{\lambda \over 2} t} f(t) = g(t) \le M_t, \qquad t\ge 0,$$
by (A.5) and the definition of $g$, we have
$$ f(t) \le e^{-{\lambda \over 2} t} \int_0^t e^{{\lambda \over 2} \tau} h(\tau)
d\tau,\qquad t\ge 0$$
and then
$$ \eqalign{\int_0^T f(t) dt &\le \int_0^T dt\int_0^t e^{-{\lambda \over 2} t} 
e^{{\lambda \over 2} \tau} h(\tau) d\tau  =\cr
&= \int_0^T d\tau \int_{\tau}^T
e^{-{\lambda \over 2} t}  e^{{\lambda
\over 2} \tau} h(\tau)  dt = \cr
&= \int_0^T e^{{\lambda \over 2} \tau} h(\tau) \Bigg(\int_{\tau}^T e^{-{\lambda
\over 2} t} dt\Bigg) d\tau=\cr
&= \int_0^T e^{{\lambda \over 2} \tau} h(\tau)
\cdot {2\over \lambda} \cdot \big( e^{-{\lambda \over 2} \tau} - e^{-{\lambda
\over 2}T}\big)d\tau=\cr
&= \int_0^T  h(\tau)
\cdot {2\over \lambda} \cdot \big( 1 - e^{-{\lambda \over
2}(T-\tau)}\big)d\tau \le  {2\over \lambda}\cdot \int_0^T  h(\tau)d\tau.}$$
So the proof is concluded.\bull 
\bigskip

\bigskip

\centerline {\bf Acknowledgments} 

The author would like to thank Prof. Alberto Bressan for suggesting the
problem and for many useful discussions.
\bigskip

\centerline {\bf References } \bigskip

\medskip
\newitem {\hbox to .8truecm {[1]} \hfill }{\maiuscoletto S.~Bianchini,
A.~Bressan,}  {\it A case study in vanishing viscosity,}
 Discr. Cont. Dyn. Sys. {\bf 7} (2001), 449-476. 

\newitem {\hbox to .8truecm {[2]} \hfill }{\maiuscoletto A.~Friedman,}
 {\it Partial differential equations of parabolic type ,}
 Prentice-Hall,  Englewood (1964). 

\newitem {\hbox to .8truecm {[3]} \hfill }{\maiuscoletto O. A. Ladyzenskaja,
 V. A. Solonnikov, N. N. Ural'ceva,}
 {\it Linear and Quasilinear Equations of Parabolic Type,}
 Translations of Mathematical Monographs, vol. 23, American Mathematical
Society (1968).

\newitem {\hbox to .8truecm {[4]} \hfill }{\maiuscoletto A. Lunardi,}
 {\it Analytic semigroups and optimal regularity in parabolic problems,}
 Birkhauser  (1995).  
\endnewitem

\end